# A non-vanishing result for the tautological ring of $\mathcal{M}_g$.

## Carel Faber

Looijenga recently proved that the tautological ring of $\mathcal{M}_g$ vanishes in degree $d > g-2$ and is at most one-dimensional in degree $g-2$, generated by the class of the hyperelliptic locus. Here we show that $\kappa_{g-2}$ is non-zero on $\mathcal{M}_g$. The proof uses the Witten conjecture, proven by Kontsevich. With similar methods, we expect to be able to prove some (possibly all) of the identities in degree $g-2$ in the tautological ring that are part of the author's conjectural explicit description of the ring.

**Lemma 1.** *The classes $\lambda_{g-1}\lambda_g$ and $ch_{2g-1}(\mathbb{E})$ vanish on the boundary $\overline{\mathcal{M}}_g - \mathcal{M}_g$.*

*Proof.* On $\Delta_0$, the Hodge bundle maps onto a trivial line bundle via the residue (after going over to a finite cover). Hence $\lambda_g = c_g(\mathbb{E})$ vanishes on $\Delta_0$. On $\Delta_i$, with $i > 0$, the Hodge bundle is the direct sum of the Hodge bundles in genus $i$ and $g-i$. Then $\lambda_g = pr_i^*\lambda_i \cdot pr_{g-i}^*\lambda_{g-i}$ and $\lambda_{g-1} = pr_i^*\lambda_i \cdot pr_{g-i}^*\lambda_{g-i-1} + pr_i^*\lambda_{i-1} \cdot pr_{g-i}^*\lambda_{g-i}$. Since in genus $h$ the relation $\lambda_h^2 = 0$ holds (an immediate consequence of the identity

$$(1 + \lambda_1 + \lambda_2 + \lambda_3 + \ldots + \lambda_h)(1 - \lambda_1 + \lambda_2 - \lambda_3 + \ldots + (-1)^h \lambda_h) = 1 \qquad (1)$$

established by Mumford), we conclude the vanishing of $\lambda_{g-1}\lambda_g$ on $\Delta_i$.

As to $ch_{2g-1}(\mathbb{E})$, the argument is even simpler: the Chern character is additive in exact sequences; but, in genus $h$, *all* the components of the Chern character of $\mathbb{E}$ in degree at least $2h$ vanish, as follows easily from the vanishing of the *even* components $ch_{2k}(\mathbb{E})$ for all $k > 0$. The vanishing of the even components is an equivalent formulation (the one that Mumford gives) of the identity (1) above.

In fact,

$$\lambda_{g-1}\lambda_g = (-1)^{g-1}(2g-1)! \cdot ch_{2g-1}(\mathbb{E}) \qquad (2)$$

(in genus $g$). We prove this below; from the formulas given, one easily proves the claimed equivalence of the two formulations of (1).

As is well-known,

$$\sum_{i=0}^{\infty} c_i t^i = \exp\left(\sum_{j=1}^{\infty} (-1)^{j+1}(j-1)! \cdot ch_j t^j\right) \ .$$

Applying this to the Hodge bundle, we find

$$\sum_{k=0}^{\infty}(2k+1)! \cdot ch_{2k+1}(\mathbb{E}) t^{2k} = \frac{d}{dt}\log\left(\sum_{i=0}^{g}\lambda_i t^i\right)$$

$$= \left(\sum_{i=0}^{g} i\lambda_i t^{i-1}\right)\left(\sum_{i=0}^{g}\lambda_i t^i\right)^{-1}$$

$$= \left(\sum_{i=0}^{g} i\lambda_i t^{i-1}\right)\left(\sum_{i=0}^{g}\lambda_i(-t)^i\right)$$

which implies
$$(2g-1)! \cdot ch_{2g-1}(\mathbb{E}) = (-1)^g(g-1)\lambda_{g-1}\lambda_g + (-1)^{g-1}g\lambda_{g-1}\lambda_g = (-1)^{g-1}\lambda_{g-1}\lambda_g$$
as well as the vanishing of the components of the Chern character of $\mathbb{E}$ in degree at least $2g$. This finishes the proof of the lemma and of (2).

**Lemma 2.** On $\overline{\mathcal{M}}_g$ :
$$\kappa_{g-2}\lambda_{g-1}\lambda_g = \frac{|b_{2g}|(g-1)!}{2^g(2g)!} \quad .$$

Here $b_k$ is the $k$-th Bernoulli number ($b_2 = \frac{1}{6}$, $b_3 = 0$, $b_4 = -\frac{1}{30}$, ...). Together, Lemma 1 and Lemma 2 imply the non-vanishing of $\kappa_{g-2}$ on $\mathcal{M}_g$ (since $b_{2g}$ doesn't vanish).

*Proof of Lemma 2.* Using Mumford's expression for the Chern character of the Hodge bundle, we find that Lemma 2 is equivalent to the identity

$$\kappa_{g-2}\left[\kappa_{2g-1} + \frac{1}{2}\sum_{h=0}^{g-1} i_{h,*}(K_1^{2g-2} - K_1^{2g-3}\cdot K_2 + \ldots + K_2^{2g-2})\right] = \frac{g!}{2^{g-1}(2g)!} \quad . \quad (3)$$

In terms of Witten's tau-classes, the left-hand-side of (3) can be rewritten as follows:

$$\langle \tau_{g-1}\tau_{2g}\rangle - \langle \tau_{3g-2}\rangle + \frac{1}{2}\sum_{j=0}^{2g-2}(-1)^j\langle \tau_{2g-2-j}\tau_j\tau_{g-1}\rangle$$
$$+ \frac{1}{2}\sum_{h=1}^{g-1}\left((-1)^{g-h}\langle \tau_{3h-g}\tau_{g-1}\rangle\langle \tau_{3(g-h)-2}\rangle + (-1)^h\langle \tau_{3h-2}\rangle\langle \tau_{3(g-h)-g}\tau_{g-1}\rangle\right) .$$

Here the first two terms come from $\kappa_{g-2}\kappa_{2g-1}$ and the first sum comes from $\Delta_0$. The second sum comes from the $\Delta_i$ with $i > 0$ ; it equals

$$\sum_{h=1}^{g-1}\frac{(-1)^{g-h}}{24^{g-h}(g-h)!}\langle \tau_{3h-g}\tau_{g-1}\rangle$$

since $\langle \tau_{3k-2}\rangle = 1/(24^k k!)$ (an immediate consequence of the Witten conjecture).

To prove the lemma, it suffices then to prove the two identities

$$\sum_{h=1}^{g}\frac{(-1)^{g-h}}{24^{g-h}(g-h)!}\langle \tau_{3h-g}\tau_{g-1}\rangle = \frac{1}{24^g g!} \quad (4)$$

and
$$\sum_{j=0}^{2g-2}(-1)^j\langle \tau_{2g-2-j}\tau_j\tau_{g-1}\rangle = \frac{g!}{2^{g-2}(2g)!} \quad . \quad (5)$$

Equation (4) follows easily from a formula I learned from Dijkgraaf ("Some facts about tautological classes", private communication (1993)):

$$\langle \tau_0 \tau(w) \tau(z) \rangle = \exp\left(\frac{(w^3 + z^3)\hbar}{24}\right) \sum_{n \geq 0} \frac{n!}{(2n+1)!} \left[\tfrac{1}{2} wz(w+z)\hbar\right]^n . \tag{6}$$

Here
$$\tau(z) = \sum_{n \geq 0} \tau_n z^n \quad \text{and} \quad \langle \mathcal{O} \rangle = \sum_{g \geq 0} \langle \mathcal{O} \rangle_g \hbar^g ,$$

so $\hbar$ is just a counting parameter, keeping track of the genus.

Indeed, for all $k \geq 1$,

$$\sum_{h=0}^{g} \frac{(-1)^{g-h}}{24^{g-h}(g-h)!} \langle \tau_0 \tau_{3h-g+k} \tau_{g-k} \rangle = 0 ,$$

since this is the coefficient of $w^{2g+k} z^{g-k} \hbar^g$ in

$$\langle \tau_0 \tau(w) \tau(z) \rangle \cdot \langle \tau_0 \tau(-w) \tau_0 \rangle = \exp\left(\frac{z^3 \hbar}{24}\right) \sum_{n \geq 0} \frac{n!}{(2n+1)!} \left[\tfrac{1}{2} wz(w+z)\hbar\right]^n ,$$

which in genus $g$ has only powers of $z$ with exponent $\geq 3i + (g-i) = g + 2i \geq g$.

Since $\langle \tau_0 \tau_a \tau_b \rangle = \langle \tau_{a-1} \tau_b \rangle + \langle \tau_a \tau_{b-1} \rangle$ we get

$$\sum_{h=0}^{g} \frac{(-1)^{g-h}}{24^{g-h}(g-h)!} \langle \tau_{3h-g} \tau_{g-1} \rangle = - \sum_{h=0}^{g} \frac{(-1)^{g-h}}{24^{g-h}(g-h)!} \langle \tau_{3h-g+1} \tau_{g-2} \rangle$$

$$= + \sum_{h=0}^{g} \frac{(-1)^{g-h}}{24^{g-h}(g-h)!} \langle \tau_{3h-g+2} \tau_{g-3} \rangle = \ldots = (-1)^{g-1} \sum_{h=0}^{g} \frac{(-1)^{g-h}}{24^{g-h}(g-h)!} \langle \tau_{3h-1} \tau_0 \rangle$$

$$= \sum_{h=1}^{g} \frac{1}{24^h h!} \frac{(-1)^{h+1}}{24^{g-h}(g-h)!} = \frac{1}{24^g g!} \left( \sum_{h=1}^{g} (-1)^{h+1} \binom{g}{h} \right) = \frac{1}{24^g g!} ,$$

which proves (4) for $g \geq 1$ (assuming (6)).

*Proof of Dijkgraaf's formula (6).* If we put $D(w, z) = \langle \tau_0 \tau(w) \tau(z) \rangle$, Witten's KdV-equation (2.33) translates into

$$\left(2w \frac{\partial}{\partial w} + 1\right) \left((w+z) D(w, z)\right) = w D(w, z) + \frac{1}{4} (w+z)^3 w D(w, z) \hbar$$
$$+ D(w, 0) z D(0, z) + 2w D(w, 0) D(0, z) .$$

It is easy to verify that the right-hand-side of (6) is the unique symmetric solution of this equation satisfying

$$D(0, z) = \langle \tau_0 \tau_0 \tau(z) \rangle = \exp\left(\frac{z^3 \hbar}{24}\right) ,$$

which finishes the proof of (6).

The proof of (5) uses three-point-functions. Writing $E(x,y,z) = \langle \tau(x)\tau(y)\tau(z) \rangle$, we find that the KdV-equation translates into

$$\left(2x\frac{\partial}{\partial x} + 1\right)\left((x+y+z)^2 E(x,y,z)\right) = E(x,0,0)(y+z)^2 E(0,y,z)$$
$$+ xE(x,y,0)zE(0,0,z) + xE(x,0,z)yE(0,y,0) + x(x+y+z)E(x,y,z)$$
$$+ 2xE(x,0,0)(y+z)E(0,y,z) + 2x(x+y)E(x,y,0)E(0,0,z)$$
$$+ 2x(x+z)E(x,0,z)E(0,y,0) + \frac{1}{4}(x+y+z)^4 xE(x,y,z)\hbar.$$

We don't need the general three-point-function, only certain coefficients of $F(w,z) := E(w,z,-z)$. The differential equation for $F(w,z)$ becomes

$$4w^2 F(w,z) + 2w^3 \frac{\partial F}{\partial w}(w,z) - \frac{1}{4}w^5 F(w,z)\hbar$$
$$= w(2w+z)D(w,z)D(0,-z) + w(2w-z)D(w,-z)D(0,z) \quad ;$$

it is clear that it has a unique solution. One verifies easily that the solution is

$$F(w,z) = \exp\left(\frac{w^3 \hbar}{24}\right) \sum_{a,b \geq 0} (w^3)^a (wz^2)^b \hbar^{a+b} \frac{(a+b)!}{2^{a+b-1}(2a+2b+2)!} \binom{a+b+1}{2a+1} \quad .$$

The coefficient of $w^g z^{2g} \hbar^g$ equals

$$\frac{(g+1)!}{2^{g-1}(2g+2)!} \quad .$$

This finishes the proof of (5) and of Lemma 2.